\documentclass[reqno]{amsart}

\usepackage{amsmath, amsfonts, amsthm, amssymb, color,  graphicx, mathrsfs, cite}
\usepackage{comment,stmaryrd}

\usepackage{ifpdf}
\ifpdf \usepackage[colorlinks=true, citecolor=blue, linkcolor=blue, urlcolor=blue]{hyperref} \fi

\theoremstyle{plain}
\newtheorem{theorem}{Theorem}[section]
\newtheorem{lemma}{Lemma}[section]
\newtheorem{proposition}{Proposition}[section]
\newtheorem{corollary}{Corollary}[section]

\theoremstyle{definition}
\newtheorem{definition}{Definition}[section]

\theoremstyle{remark}

\newtheorem{example}{Example}[section]

\numberwithin{equation}{section}
\allowdisplaybreaks

\usepackage{ifpdf}
\ifpdf \usepackage[colorlinks=true, citecolor=blue, linkcolor=blue, urlcolor=blue]{hyperref} \fi

\def\thrm{\begin{theorem}}
\def\thrml#1{\begin{theorem}\label{#1}}
\def\ethrm{\end{theorem}}
\def\lmm{\begin{lemma}}
\def\lmml#1{\begin{lemma}\label{#1}}
\def\elmm{\end{lemma}}
\def\dfntn{\begin{definition}}
\def\dfntnl#1{\begin{definition}\label{#1}}
\def\edfntn{\end{definition}}
\def\crllr{\begin{corollary}}
\def\crllrl#1{\begin{corollary}\label{#1}}
\def\ecrllr{\end{corollary}}
\def\xmpl{\begin{example}}
\def\xmpll#1{\begin{example}\label{#1}}
\def\exmpl{\end{example}}
\def\nmrt{\begin{enumerate}}
\def\enmrt{\end{enumerate}}
\def\qtn{\begin{equation}}
\def\qtnl#1{\begin{equation}\label{#1}}
\def\eqtn{\end{equation}}
\def\prpstn{\begin{proposition}}
\def\prpstnl#1{\begin{proposition}\label{#1}}
\def\eprpstn{\end{proposition}}

\def\tm#1{\item[{\rm (#1)}]}
\def\proof{{\bf Proof}.\ }
\def\eprf{\hfill$\square$}

\DeclareMathOperator{\aut}{Aut}

\DeclareMathOperator{\Aaut}{Aut_{alg}}
\DeclareMathOperator{\AGL}{AGL}

\DeclareMathOperator{\inv}{Inv}

\newcommand{\FT}{\mathbf{T}}

\def\fX{{\frak X}}

\newcommand{\Om}{\Omega}

\def\qaq{\quad\text{and}\quad}

\def\mF{{\mathbb F}}

\def\Tr{{\rm Tr}}
\def\cX{\mathcal {X}}
\def\cY{\mathcal {Y}}

\def\qaq{\quad\text{and}\quad}

\DeclareMathOperator{\AGaL}{A{\rm \Gamma}L}
\DeclareMathOperator{\diag}{Diag}
\DeclareMathOperator{\PGL}{PGL}
\DeclareMathOperator{\GL}{GL}

\DeclareMathOperator{\orb}{Orb}

\DeclareMathOperator{\PGaL}{P{\rm \Gamma}L}
\DeclareMathOperator{\PSL}{PSL}

\DeclareMathOperator{\sym}{Sym}
\DeclareMathOperator{\WL}{WL}

\begin{document}

\title[A characterization of exceptional pseudocyclic schemes]{A characterization of exceptional pseudocyclic association schemes by multidimensional intersection numbers}

\author{Gang Chen}
\address{School of Mathematics and Statistics, Central China Normal University, Wuhan, China}
\email{chengangmath@mail.ccnu.edu.cn}

\author{Jiawei He}
\address{School of Mathematics and Statistics, Central China Normal University, Wuhan, China}
\email{hjwywh@mails.ccnu.edu.cn}

\author{Ilia Ponomarenko}
\address{School of Mathematics and Statistics of Central China Normal University, Wuhan, China and Steklov Institute of Mathematics at St. Petersburg, Russia}
\email{inp@pdmi.ras.ru}

\author{Andrey Vasil'ev}
\address{Sobolev Institute of Mathematics, Novosibirsk, Russia}
\email{vasand@math.nsc.ru}

\thanks{The first two authors is supported by the NSFC grant No. 11971189. The third and fourth authors are supported by the RFBR grant No. 18-01-00752}

\begin{abstract}
Recent classification of $\frac{3}{2}$-transitive permutation groups leaves us with three infinite families of groups which are neither $2$-transitive, nor Frobenius, nor one-dimensional affine. The groups of the first two families correspond to special actions of $\PSL(2,q)$ and $\PGaL(2,q)$, whereas those of the third family are the affine solvable subgroups of $\AGL(2,q)$ found by D.~Passman in~1967. The association schemes of the groups in each of these families are known to be pseudocyclic. It is proved that apart  from three particular cases, each of these exceptional pseudocyclic schemes is characterized up to isomorphism by the tensor of its $3$-dimensional intersection numbers.
\medskip

\end{abstract}


\date{}

\maketitle


\section{Introduction}

In the late 1960s, H.~Wielandt proposed a method for studying permutation groups via invariant relations. Later, D.~Higman  axiomatized a part of this method (connected with binary relations) by introducing a new object called  a coherent configuration~\cite{Hig}. The coherent configuration of a permutation group $G$ is formed by the orbits of  the induced action of~$G$ on the Cartesian square of the underlying set of points (for exact definitions, see Section~\ref{sect:pre}). Looking only at the parameters of this coherent configuration, the so-called intersection numbers, one can easily determine whether the original group is transitive,  primitive, $2$-transitive, etc. For example, the transitivity of a group $G$ means exactly that the coherent configuration of $G$ is an association scheme.\medskip

The concept of  {\it pseudocyclic} (association) scheme goes back to research of D.~Mesner~\cite{M64}, related with constructing of designs and strongly regular graphs; the defining property of such a scheme is that the ratio of the multiplicity and degree of its nonprincipal irreducible character does not depend on the choice of the character. It was proved in~\cite[Theorem~3.2]{MP1} that this is equivalent to a certain relation for intersection numbers, see Subsection~\ref{240720a}.\medskip

The class of pseudocyclic schemes contains all Frobenius schemes, i.e., the coherent configurations of the Frobenius groups, and, moreover,  every  pseudocyclic scheme of rank sufficiently large comparing with its degree is Frobenius~\cite[Theorem~1.1]{MP1}.  Thus the pseudocyclic schemes can be considered as combinatorial analogs of the Frobenius groups. It should be mentioned that the analogy is not complete, because there exist schurian (i.e., those associated with permutation groups) pseudocyclic schemes which are not Frobenius, as well as non-schurian pseudocyclic schemes~\cite[Example 2.6.15]{CP1}. \medskip

It is well known that every Frobenius group is $\frac{3}{2}$-transitive, i.e., is transitive  and all the orbits of the stabilizer of a point $\alpha$, other than~$\{\alpha\}$,  have the same size greater than~$1$. It immediately follows that so is the automorphism group of any Frobenius scheme. Moreover, the above mentioned relation for intersection numbers implies that the automorphism group of a schurian pseudocyclic scheme is  also $\frac{3}{2}$-transitive. A recent classification of  $\frac{3}{2}$-transitive permutation groups 
shows that in most cases the coherent configuration of a $\frac{3}{2}$-transitive group is pseudocyclic, see Subsection~\ref{110620a}. We cite a part of this classification in the following theorem, see~\cite[Corollaries 2,3]{LPS}.

\thrml{140520c}
Let $G$ be a $\frac{3}{2}$-transitive permutation group of degree~$n$. Assume that neither $G$ is $2$-transitive or  Frobenius nor $G\le\AGaL(1,q)$ for some~$q$. Then apart from finitely many cases, 
\nmrt
\tm{1} $n=q(q-1)/2$,  $q=2^d\ge 8$,  and either $G=\PSL(2,q)$, or $d$ is prime and $G=\PGaL(2,q)$,
\tm{2} $n=q^2$,  $q$ is odd, and  $G\le \AGL(2,q)$ is the affine group with point stabilizer of order $4(q-1)$, consisting of all monomial matrices of determinant $\pm 1$.
\enmrt
\ethrm

The association schemes of the groups in statement~(1) of  Theorem~\ref{140520c} appeared in Master Thesis of H.~Hollmann (1982); such a scheme is called a {\it large} or {\it small  Hollmann scheme} depending on whether $G=\PSL(2,q)$ or $G=\PGaL(2,q)$.\footnote{According to the Galois correspondence between permutation groups and coherent configurations~\cite[Section~2.2]{CP1}, the smaller groups correspond to larger coherent configurations}  These schemes have been studied in~\cite{HX2}; in particular, it was proved there that both are pseudocyclic.  The group in statement~(2) of  Theorem~\ref{140520c} appeared  in D.~Passman's characterization of solvable $\frac{3}{2}$-transitive groups~\cite{P67}. The association scheme of this group, the {\it Passman scheme}, is also pseudocyclic~\cite{MP1}.\medskip

The goal of the present paper is to establish combinatorial characterizations of the Hollmann and Passman schemes; from the point of view of Theorem~\ref{140520c},  they can naturally be considered as exceptional. One of the best possible combinatorial characterization of an association scheme is obtained when the scheme in question is determined up to (combinatorial) isomorphism by its intersection numbers; in this case the scheme is called separable. However, most of association  schemes are not separable. In \cite{EP4}, multidimensional intersection numbers and separability number~$s(\cX)$ of a coherent configuration $\cX$ have been introduced and studied (see also~\cite[Section~3.5 and~4.2]{CP1}). According to the definition, $s(\cX)\le m$ if and only if $\cX$  is determined up to isomorphism by its $m$-dimensional intersection numbers; thus $s(\cX)=1$ if and only if $\cX$ is separable. It was proved in~ \cite{EP4} that $s(\cX)=1$ or~$2$ if $\cX$ is the scheme of a classical distance regular graph; later, the estimate $s(\cX)\le 3$ has been established for any cyclotomic scheme $\cX$ over finite field.

\thrml{130520a}
Let $\cX$ be a large Hollmann scheme. Then  $s(\cX)\le 2$.
\ethrm

The proof of Theorem~\ref{130520a} is given in Section~\ref{Hs}. The difficult step in the proof is to verify that the one point extension of the scheme $\cX$ (which is a combinatorial analog of a one point stabilizer of permutation group) is a coherent configuration of the stabilizer of this point in $\aut(\cX)$. In proving this fact we use the formulas for the intersection numbers of $\cX$, which were calculated in~\cite{HX2}.\medskip

The proof of the following two theorems is based on Theorem~\ref{201444a} (see Section~\ref{lb}), giving a sufficient condition for an arbitrary coherent configuration $\cX$ to be {\it partly regular}, i.e., to be the coherent configuration of a permutation group having a faithful regular orbit. Using this sufficient condition we are able to show that if~$\cX$ is the small Hollmann scheme (apart from several exceptions) or the Passman scheme, then a two point extension of $\cX$ is partly regular. Modulo known results, this immediately implies that $s(\cX)\le 3$.

\thrml{130520x}
Let $\cX$ be a small Hollmann scheme of degree $q(q-1)/2$, where $q=2^d$ with prime $d\ne 7,11,13$. Then $s(\cX)\le 3$.
\ethrm

In the three exceptional cases of Theorem~\ref{130520x}, the sufficient condition given in Theorem~\ref{201444a} does not work. It seems that the conclusion of Theorem~\ref{130520x} is also true for them. However, the corresponding schemes are too large to check this  statement by a direct calculation.

\thrml{010620x}
Let $\cX$ be a Passman scheme. Then $s(\cX)\le 3$. 
\ethrm

Throughout the paper, we actively use the notation, concepts, and statements from the theory of coherent configurations. All of them can be found in the monograph~\cite{CP1}. In Section~\ref{sect:pre}, we give a brief extract from the theory of coherent configurations that is relevant for this paper.

\medskip

{\bf Notation.}

For a prime power $q$, $\mF_q$ is a finite field of order~$q$.

Throughout the paper, $\Omega$ is a finite set. 

The diagonal of the Cartesian product $\Omega\times\Omega$ is denoted by $1_\Omega$;  for $\alpha\in \Omega$, we set~$1_\alpha:=1_{\{\alpha\}}$. 

For $r\subseteq \Omega\times\Omega$, we set $r^*=\{(\beta,\alpha) :\ (\alpha,\beta)\in r\}$ and $\alpha r=\{\beta\in \Omega:\ (\alpha,\beta)\in r\}$, $\alpha\in \Omega$.

For relations $r,s\subseteq \Omega\times\Omega$, we set $r\cdot s=\{(\alpha,\beta):\ (\alpha,\gamma)\in r,\  (\gamma,\beta)\in s\}$. %

For a set $S$ of relations on $\Omega$, we define $S^*=\{s^*:\ s\in S\}$ and put $S^\cup$ to be the set of all unions of the relations of $S$. 

\section{Coherent configurations}\label{sect:pre}

\subsection{Rainbows}

Let $\Omega$ be a finite set and  $S$ a partition of $\Omega\times\Omega$. A pair $\cX=(\Omega,S)$ is called
a {\it rainbow} on $\Omega$ if
$$
1_\Omega\in S^\cup, \textrm{~and~} S^*=S.
$$
The elements of the sets $\Omega$, $S=S(\cX)$, and $S^{\cup}$ are called, respectively, the {\it points},  {\it basis relations},  and  {\it relations} of~$\cX$. The numbers $|\Omega|$ and $|S|$ are called the {\it degree} and {\it rank} of~$\cX$, respectively. The unique basic relation containing a pair $(\alpha,\beta)\in\Omega\times\Omega$ is denoted by $r_\cX(\alpha,\beta)$; we omit the subscript~$\cX$ wherever it does not lead to misunderstanding.\medskip

A set $\Delta\subseteq\Omega$ is called a  {\em fiber} of a rainbow~$\cX$ if $1_\Delta\in S$; the set of all fibers is denoted by $F:=F(\cX)$. The point set~$\Omega$ is the disjoint union of fibers. If $\Delta$ is a union of fibers, then the pair
\[
\cX_\Delta=(\Delta,S_\Delta)
\]
is a rainbow, where $S_\Delta$ consists of all $s_\Delta=s\cap(\Delta\times\Delta)$, $s\in S$.\medskip

Let $\cX=(\Omega,S)$ and  $\cX'=(\Omega',S')$ be rainbows. A bijection  $f:\Omega\to\Omega'$ is called a {\it combinatorial isomorphism} (or simply isomorphism) from $\cX$ to $\cX'$ if $S^f=S'$. When $\cX=\cX'$, the set of all these isomorphisms form a permutation group on~$\Omega$. This group has a (normal) subgroup
$$
\aut(\cX)=\{f\in\sym(\Omega)\!:\ s^f=s\ \,\text{for all}\ \, s\in S\}, 
$$
called the  {\it automorphism group} of~$\cX$.

\subsection{Coherent configurations}
A rainbow $\cX=(\Omega,S)$ is called a  {\it coherent configuration} if for any $r,s,t\in S$, the number
$$
c_{rs}^t=|\alpha r\cap \beta s^*|
$$
does not depend on the choice of $(\alpha,\beta)\in t$;  the numbers~$c_{rs}^t$ are called the {\it intersection numbers} of~$\cX$. If, in addition, $1_\Omega\in S$, then the coherent configuration $\cX$ is said to be {\it homogeneous}, an {\it association scheme}, or just  a {\it scheme}. A scheme $\cX$ is called {\it symmetric} if $s=s^*$ for all $s\in S$.\medskip

Let $\cX$ be a coherent configuration. Then for any $s\in S$, there exist uniquely determined $\Delta,\Gamma\in F$ such that $s\subseteq \Delta\times\Gamma$. Denote by $S_{\Delta,\Gamma}$ the set of all~$s$ contained in $\Delta\times\Gamma$. Then the union
$$
S=\bigcup_{\Delta,\Gamma\in F}S_{\Delta,\Gamma}
$$
is disjoint. The positive integer $|\delta s|$, $\delta\in\Delta$, equals the intersection number~$c_{ss^*}^{1_\Delta}$,
and hence does not depend on the choice of~$\delta$. It is called the  {\it valency} of $s$ and denoted by $n_s$.  In homogeneous case, $n_s=n_{s^*}$ and also  
\qtnl{150309c}
n_tc_{rs}^{t^*}=n_rc_{st}^{r^*}=n_sc_{tr}^{s^*},\qquad r,s,t\in S.
\eqtn

A basis relation  $s\in S$ is called a  {\it matching} if  $n_{s^{}}=n_{s^*}=1$.  Note that a matching $s\in S_{\Delta,\Gamma}$ defines a bijection from $\Delta$ to  $\Gamma$, taking $\delta\in\Delta$ to the unique point of the singleton $\delta s$. Furthermore, one can see that if $r\in S$ and  $t=s\cdot r$ (respectively, $t=r\cdot s$)  is nonempty, then $t\in S$.\medskip

Let $G$ be a permutation group on~$\Omega$. Denote by $(\alpha,\beta)^G$ the orbit of the induced action of~$G$ on~$\Omega\times\Omega$, that contains the pair~$(\alpha,\beta)$. Then 
$$
\inv(G)=\inv(G,\Omega)=(\Omega,\{(\alpha,\beta)^G:\ \alpha,\beta\in\Omega\})
$$
is a coherent configuration; we say that $\inv(G)$ is the coherent configuration associated with $G$.  A coherent configuration $\cX$ is said to be {\it schurian} if $\cX=\inv(\aut(\cX))$. 

\subsection{Separability}
Let $\cX=(\Omega,S)$ and $\cX'=(\Omega',S')$ be coherent configurations. A bijection $\varphi:S\to S',\ s\mapsto s'$, is called an {\it algebraic isomorphism} from~$\cX$ to~$\cX'$ if
\qtnl{f041103p1}
c_{rs}^t=c_{r's'}^{t'},\qquad r,s,t\in S.
\eqtn
When $\cX=\cX'$, the set of all such $\varphi$ forms a subgroup of $\sym(S)$, denoted by $\Aaut(\cX)$.\medskip

Each isomorphism~$f$ from~$\cX$ to~$\cX'$ induces an algebraic isomorphism from~$\cX$ to~$\cX'$, which maps $r\in S$ to $r^f\in S'$. A coherent configuration~$\cX$ is said to be  {\it separable} if  every algebraic isomorphism from $\cX$ is induced by a suitable bijection (which is an isomorphism of the coherent configurations in question).\medskip

The algebraic isomorphism $\varphi$ induces a bijection from $S^\cup$ to $(S')^\cup$: the union $r\cup s\cup\cdots$ of basis relations of $\cX$ is mapped to $r'\cup s'\cup\cdots$. This bijection is also denoted by $\varphi$. It preserves the dot product, i.e., $\varphi(r\cdot s)=\varphi(r)\cdot\varphi(s)$ for all $r,s\in S$.\medskip

\subsection{Coherent closure}

There is a natural partial order\, $\le$\, on the set of all rainbows on the same set~$\Omega$. Namely, given two such rainbows  $\cX$ and $\cX'$, we set
$$
\cX\le\cX'\ \Leftrightarrow\ S(\cX)^\cup\subseteq S(\cX')^\cup.
$$
The minimal and maximal elements with respect to this order are the {\it trivial} and {\it discrete} coherent
configurations, respectively: the basis relations of the former are the reflexive relation $1_\Omega$ and (if $|\Omega|>1$) its complement in $\Omega\times\Omega$, whereas the basis relations of the latter are singletons.\medskip

The functors $\cX\to\aut(\cX)$  and $G\to\inv(G)$ form a Galois correspondence between the posets of coherent configurations and permutation groups on the same set, i.e.,  
$$
\cY\le \cX \Rightarrow  \aut(\cY)\ge \inv(\cX)\qaq L\le K\Rightarrow \inv(L)\ge \inv(K),
$$
and 
$$
\aut(\inv(\aut(\cX)))=\aut(\cX)\qaq \inv(\aut(\inv(G)))=\inv(G).
$$

The {\it coherent closure} $\WL(T)$ of a set $T$ of relations on $\Omega$, is defined to be the smallest coherent configuration
on $\Omega$, for which $T$ is a set of relations. The {\it point extension} $\cX_{\alpha,\beta,\ldots}$ of the rainbow~$\cX$ with respect to the points $\alpha,\beta,\ldots\,\in\Omega$ is defined to be $\WL(T)$, where $T$ consists of $S(\cX)$ and the relations $1_\alpha,1_\beta,\ldots$. In other words, $\cX_{\alpha,\beta,\ldots}$ is the smallest coherent configuration on $\Omega$ that is larger than or equal to~$\cX$ and has singletons $\{\alpha\},\{\beta\},\ldots$ as fibers.

\subsection{Multidimensional intersection numbers}

The theory of multidimensional extensions of coherent configurations has been developed in \cite{EP4}, see also~\cite[Section~3.5]{CP1}. \medskip

Let $m\ge 1$ be an integer. The $m$-extension of a coherent configuration $\cX$ on~$\Omega$ is defined to be the smallest coherent configuration on $\Omega^m$, which contains the Cartesian $m$-power of $\cX$ and for which the set $\diag(\Omega^m)$ is the union of fibers. The intersection numbers of the $m$-extension are called the {\it $m$-dimensional intersection numbers} of the configuration~$\cX$.  If  $m=1$, then the $m$-extension of~$\cX$ coincides with~$\cX$ and the  $m$-dimensional intersection numbers of $\cX$ are the ordinary intersection numbers.\medskip

An algebraic isomorphism $\varphi$ from $\cX$ to $\cX'$ is said to be {\it $m$-dimensional} if it can be extended to an algebraic isomorphism from the $m$-extension of~$\cX$ to that of~$\cX'$, that takes $\diag(\Omega^m)$ to $\diag({\Omega'}^m)$. The {\it separability number} $s(\cX)$ of the coherent configuration $\cX$ is defined to be the  smallest positive integer $m$ for which every algebraic $m$-dimensional isomorphism from $\cX$ is induced by some isomorphism. Thus, the equality $s(\cX)=m$ expresses the fact that $\cX$ is determined up to isomorphism by  its tensor of the $m$-dimensional intersection numbers. The following statement was proved in \cite[Theorem~4.6(1)]{EP4}.

\lmml{030620d}
Let $\cX$ be a coherent configuration. Then $s(\cX)\le s(\cX_\alpha)+1$ for any point $\alpha$ of~$\cX$.
\elmm

\subsection{Pseudocyclic schemes}\label{240720a}
Let $\cX=(\Omega,S)$ be a  coherent configuration. The {\it indistinguishing number} of a relation $s\in S(\cX)$ is defined to be the sum  $c(s)$ of the intersection numbers $c_{rr^*}^s$, $r\in S$. For each pair $(\alpha,\beta)\in s$,  we have $c(s)=|c(\alpha,\beta)|$, where
\qtnl{100814c}
c(\alpha,\beta)=\{\gamma\in\Omega:\ r(\gamma,\alpha)=r(\gamma,\beta)\}.
\eqtn
The maximum $c(\cX)$ of the numbers $c(s)$, where $s$ runs over the set of all irreflexive  basis relations of~$\cX$, is called the {\it indistinguishing number} of~$\cX$. It is easily seen that $c=0$ if and only if $n_s=1$ for each $s\in S$.\medskip

Assume that $\cX$ is a scheme. In accordance with  \cite[Theorem 3.2]{MP1}, $\cX$ is  {\it pseudocyclic} of valency $k$ if the equalities  
$$
c(s)+1=k=n_s
$$
hold for all irreflexive  $s\in S$.  The class of pseudocyclic schemes includes all (homogeneous) coherent configurations associated with regular or Frobenius  groups.

\subsection{Partly regular coherent configurations}

A coherent configuration $\cX$ is said to be {\it partly regular} if there exists a point $\alpha \in \Omega$ such that $|\alpha s|\le 1$ for all $s \in S$; the point $\alpha$ is said to be~{\it regular}.  When all the points of $\cX$ are regular, we say that $\cX$ is {\it semiregular}, and {\it regular} if $\cX$ is a scheme. Thus, $\cX$ is semiregular if and only if $c(\cX)=0$.\medskip
 
The following statement taken from \cite[Theorem~3.3.19]{CP1}  shows, in particular, that the partly regular (respectively, semiregular, regular) coherent configurations are in one-to-one correspondence with those of the form $\inv(G)$, where $G$ is a permutation group having a faithful orbit  (respectively, $G$ is semiregular, regular).
 
\thrml{200520i}
Every partly regular coherent configuration $\cX$ is schurian and separable. In particular, $s(\cX)=1$.
\ethrm

The key point in the proof of Theorem~\ref{200520i} is the lemma below \cite[Lemma~3.3.20]{CP1}; it is also used in the proof of Theorem~\ref{130520a}.

\lmml{180520i}
Let $\cX$ be a coherent configuration and $\Delta$ a union of fibers of~$\cX$. Assume that for every $\Gamma\in F$ there exists $s\in S_{\Delta,\Gamma}$ such that $n_s =1$. Then
\nmrt
\tm{1} the restriction mapping $\aut(\cX) \rightarrow \aut(\cX_\Delta)$  is a group isomorphism,
\tm{2} $\cX$ is schurian and separable whenever $\cX_\Delta$ is schurian and separable.
\enmrt
\elmm

\section{Large  Hollmann schemes}\label{Hs}

\subsection{General properties.}\label{260720x}
Throughout this section, $d\geq 3$ is an integer, $q = 2^d $, and $G=\PSL(2,q)$ the permutation group of degree~$n=q(q-1)/2$ from  Theorem~\ref{140520c}(1).  The lemma below immediately follows from Theorem~1.2(iii) and Lemma~6.2 proved in~\cite{BGLPS}.

\lmml{250720a}
For any point $\alpha$, we have $G_\alpha=D_{2(q+1)}$. Moreover,
\qtnl{160520s}
|\Delta|=q+1,\qquad \Delta\in\orb(G_\alpha),\ \Delta\ne\{\alpha\}.
\eqtn
\elmm

Let $\cX=\inv(G)$ be the large Hollmann scheme. We need to compare $\cX$ with  symmetric pseudocyclic scheme $\cX'$ of degree $n$ and valency $q+1$, associated with the group $\PGL(2,q)=\PSL(2,q)$ and studied in~\cite{HX2}.\footnote{Not only to justify the name used first in~\cite{MP1} but also to prove Proposition~\ref{160520i} below. } To this end, we note that $\aut(\cX)=G$, because as was observed in~\cite[Lemma~2.10]{VC}, the group~$G$ is $2$-closed. On the other hand, $\cX'$ is also associated with~$G$ and has the same valencies as~$\cX$ by Lemma~\ref{250720a}. Thus  by Theorem~\ref{140520c}, we have $\aut(\cX')=G=\aut(\cX)$ and hence
\qtnl{160520w}
\cX'=\inv(\aut(\cX'))=\inv(\aut(\cX))=\cX.
\eqtn

\prpstnl{160520i}
The large Hollmann scheme $\cX$ is symmetric and pseudocyclic of  degree $q(q-1)/2$, rank $q/2$, and valency~$q+1$.  Moreover,
\qtnl{250720b}
\aut(\cX)=G\qaq \aut(\cX_\alpha)=G_\alpha=D_{2(q+1)}\ \,\text{for all}\ \,\alpha.
\eqtn
\eprpstn
\proof By the remark before the proposition, we need to verify the second equality in~\eqref{250720b} only. By \cite[Proposition 3.3.3(1)]{CP1}, we have $\aut(\cX_\alpha)=\aut(\cX)_\alpha$. Thus the required statement immediately follows from  the first equality in~\eqref{250720b} and Lemma~\ref{250720a}.\eprf\medskip

Equality \eqref{160520w} allows us to use formulas for the intersection numbers of the scheme $\cX'$, given in~\cite[Theorem 2.2]{HX2}.  Namely, let $S=S(\cX)$ and 
$$
\FT_0=\{x \in \mF_{2^d}:\ \Tr(x)=0\},
$$ 
where $\Tr(x)$ is the trace of $x$  over the prime subfield of the field $\mF_{2^d}$. Then there is a bijection $\FT_0 \rightarrow S, x \mapsto s_{x}$, such that $s_0$ is reflexive and
\qtnl{250720c}
c_{s_{x}, s_{y}}^{s_{z}}=1\quad\Leftrightarrow\quad\ \Tr(xz)=0\ \,\text{and}\ \,x+y+z=0.
\eqtn
As is easily seen,  $\FT_0$ is as a linear space of dimension $d-1$ over $\mF_2$. 

\subsection{One point extension.} Let us analyze the  extension $\cX_\alpha$ of the large Hollmann scheme~$\cX$ with respect to a  point~$\alpha$. Since the scheme $\cX$ is schurian, each fiber of the coherent configuration $\cX_\alpha$ is of the form $\Delta=\alpha s$ for some $s\in S$ \cite[Theorem~3.3.7]{CP1}. When $s=s_x$ for some $x\in\FT_0$, the fiber $\Delta$ is denoted by $\Delta_x$. Thus,
$$
F(\cX_\alpha)=\{\Delta_x:\ x\in\FT_0\}. 
$$

\thrml{170520w1}
Let $x$ and $y$ be nonzero elements of $\FT_0$. Then the set $S(\cX_\alpha)_{\Delta_x,\Delta_y}$ contains a matching.
\ethrm
\proof  We need auxiliary lemmas.
 
\lmml{170520w}
Theorem~$\ref{170520w1}$ holds whenever $\Tr(xy)=0$.
\elmm
\proof Let $z=x+y$. Then obviously  $z\in \FT_0$. Moreover 
$$
\Tr(xz)=\Tr(x^2+xy)=\Tr(x^2)+\Tr(xy)=0.
$$ 
By  formula~\eqref{250720c}, this implies that $c_{s_x, s_y}^{s_z}=1$. If $z=0$, then $x=y$ and $1_{\Delta_x}$ is a desired matching. Assume that $z$ is nonzero. Then $n_{s_x}=n_{s_z}$, because $\cX$ is a pseudocyclic  scheme (Proposition~\ref{160520i}). Since $\cX$ is also symmetric, we have
$$
c_{s_y,s_z}^{s_x}= \frac{n_{s_z}}{n_{s_x}}c_{s_x,s_y}^{s_z}=1,
$$
see \eqref{150309c}. Therefore if $r=s_z\cap (\alpha s_x\times \alpha s_y)$, then $|\beta r|=1$ for all $\beta\in \Delta_x$, and $|\beta r^*|=1$ for all $\beta\in \Delta_y$ (here we use the fact  that $|\Delta_x|=|\Delta_y|$). Since $r$ is a relation of $\cX_\alpha$ \cite[Lemma~3.3.5]{CP1}, this implies that $r$ belongs to $S(\cX_\alpha)_{\Delta_x,\Delta_y}$. Thus,~$r$ is a required matching.\eprf\medskip

Let us define a graph $\fX$ with nonzero elements of $\FT_0$ as the vertices and in which two distinct vertices $x$
and $y$ are adjacent if and only if $\operatorname{Tr}(x y)=0 .$ One can see that $\fX$ is an undirected graph with exactly $|\FT_0|-1$ vertices.

\lmml{4151533a}
The graph $\fX$ is connected.
\elmm
\proof Given $x \in \FT_0$,  the mapping
$$
f_x:\FT_0\to\mF_{2},\ y \mapsto\Tr(x y)
$$
is linear. If $x\ne 0$ and $\ker(f_x)=\FT_0$,  then $x$ is adjacent in $\fX$ with all other vertices and hence $\fX$ is a connected graph. Thus without loss of generality, we may assume that
\qtnl{170520d}
|\ker(f_x)|=\frac{|\FT_0|}{2}=\frac{q}{4}
\eqtn
for all nonzero $x\in\FT_0$.  Then the vertex $x$ has exactly $k=\frac{q}{4}-2$ neighbors in $\fX$ other than~$x$ (note that 
$0$ belongs to $\FT_0$, but not a vertex of $\fX$).  Consequently, $\fX$ is a regular graph of valency $k$ and 
$$
|\FT_0|-1=\frac{q}{2}-1=2 k+3
$$ 
vertices.\medskip

Assume on the contrary that the graph $\fX$ is not connected. Each component of~$\fX$ has  at least $k+1$ vertices. Consequently, $\fX$ has exactly two components: one with $k+1$ vertices and another one with $k+2$ vertices. Let $x$ and $y$ be distinct nonadjacent vertices of the second component. Since $\fX$ is regular of valency $k$, these vertices have exactly $k$ common neighbors in $\fX$. Consequently,
\qtnl{170529u}
|\ker(f_x) \cap \ker(f_y)|=k+1=\frac{|\FT_0|}{2}-1.
\eqtn

On the other hand, $\ker(f_x) \cap \ker(f_y)$ is a linear subspace of $\FT_0$ other than $\FT_0$, see~\eqref{170520d}. If it is a hyperplane, then $\ker(f_x)=\ker(f_y)$ and hence $x$ and $y$ are adjacent, a contradiction. Therefore, the codimension of $\ker(f_x) \cap \ker(f_y)$ is at least~$2$. Thus, 
\qtnl{170529v}
|\ker(f_x) \cap \ker(f_y)|\le \frac{|\FT_0|}{4}.
\eqtn
Comparing inequalities~\eqref{170529u} and \eqref{170529v}, we obtain $4\ge |\FT_0|=q/2=2^{d-1}$. It follows that $d=3$. A straightforward computer calculation shows that in this case the graph~$\fX$ is connected.\eprf\medskip

Let us return to the proof of Theorem~\ref{170520w1}. By Lemma~\ref{4151533a}, the vertices $x$ and~$y$ of the graph~$\fX$ are connected by a path 
$$
x=x_0,x_1,\ldots,x_k=y,
$$ 
where $k\ge 1$. For $i=0,\ldots,k-1$, the vertices $x_i$ and $x_{i+1}$ are adjacent and hence $\Tr(x_ix_{i+1})=0$. Denote by $s_i$ the matching in $S(\cX_\alpha)_{\Delta_{x_i},\Delta_{x_{i+1}}}$ the existence of which is guaranteed by Lemma~\ref{170520w}. Then the dot product
$$
s=s_0\cdot s_1\cdots s_{k-1}
$$
is a desired matching belonging to $S(\cX_\alpha)_{\Delta_x,\Delta_y}$.\eprf

\crllrl{250720f}
Let $\Delta=\Delta_x$ for nonzero $x\in\FT_0$. Then the coherent configuration $\cY=(\cX_\alpha)_\Delta$ is schurian and separable. Moreover, the extension of  $\cY$ with respect to at least one point is partly regular.
\ecrllr
\proof By Proposition~\ref{160520i}, we have $\aut(\cX_\alpha)=D_{2(q+1)}$. Furthermore, the hypothesis of Lemma~\ref{180520i} is satisfied for $\cX=\cX_\alpha$ by  Theorem~\ref{170520w1}. Thus by  statement~(1) of that lemma, we have 
\qtnl{200520a}
H:=\aut(\cY)\cong\aut(\cX_\alpha)= D_{2(q+1)}. 
\eqtn
On the other hand, $|\Delta|=q+1$ by formula~\eqref{160520s}. Consequently, the group $H$ contains a normal regular cyclic subgroup $C$ of order $q+1$. In terms of \cite[Section~4.4]{CP1}, this means that $\cY$ is isomorphic to a normal circulant  scheme. The radical of such a scheme, being a subgroup of the group~$C$, is of order at most~$2$; this follows from the implication (1)$\Leftrightarrow$(3) in~\cite[Theorem~6.1] {EP1}. Since the number $|C|=q+1=2^d+1$ is odd, the radical is trivial. Thus, the scheme~$\cY$ is schurian by \cite[Corollary 4.4.3]{CP1}, and every its extension with respect to at least one point  is partly regular by~\cite[Theorem~4.4.7]{CP1}.\medskip

It remains to verify that $\cY$ is separable. Since $Y$ is schurian by above, we have $\cY=\inv(\aut(\cY))=\inv(H)$. By virtue of~\eqref{200520a}, this means that $\cY$ is the coherent configuration associated with $D_{2(q+1)}$. Thus, the required statement follows from~\cite[Exercise~2.7.33]{CP1}.\eprf

\subsection{Proof of Theorem~\ref{130520a}.} 

By Lemma~\ref{030620d}, it suffices to verify that a one point extension of a large Hollmann scheme is separable. But this immediately follows from Theorem~\ref{030620i} below.

\thrml{030620i}
Let $q=2^d$ where $d\ge 3$. Then the extension of the large Hollmann scheme~of  degree $q(q-1)/2$ with respect to at least one point is schurian and separable. 
\ethrm
\proof Let $\cX'$ be the extension of the large Hollmann scheme $\cX$ with respect to $m\ge 1$ points $\alpha=\alpha_1,\alpha_2,\ldots,\alpha_m$. Let $x\in \FT_0$ be nonzero  and $\Delta=\Delta_x$.  Then the hypothesis of Lemma~\ref{180520i} is satisfied for $\cX=\cX'$. Indeed, each $\Gamma\in F(\cX')$ other than $\{\alpha\}$ is contained in $\Delta_y$ for some nonzero $y\in \FT_0$. By Theorem~\ref{170520w1}, there is a matching $s'\in S(\cX_\alpha)_{\Delta_x,\Delta_y}$, and as the required relation~$s$ one can take $s'\cap(\Delta\times\Gamma)$.\medskip

By Lemma~\ref{180520i}(2), the coherent configuration $\cX'$ is schurian and separable whenever so is $\cX'_\Delta$. If  $m=1$, then $\cX'=\cX_\alpha$ and we are done by Corollary~\ref{250720f}. Let $m>1$. We claim that there exist $\beta_2,\ldots,\beta_m\in\Delta$ such that 
\qtnl{210520j}
\cX'=\cX_{\alpha,\beta_2,\ldots,\beta_m}.
\eqtn
Indeed, without loss of generality we may assume that none of the $\alpha_i$, $i>1$, equals~$\alpha$. By Theorem~\ref{170520w1}, there is a matching $s_i\in S(\cX_\alpha)_{\Delta_i,\Delta}$, where $\Delta_i$ is the fiber of $\cX_\alpha$, containing~$\alpha_i$.  Then $\alpha_i s_i=\{\beta_i\}$ for some $\beta_i\in \Delta$.  It follows that for any extension of $\cX_\alpha$, each or none of the two singletons $\{\alpha_i\}$ and~$\{\beta_i\}$ is a fiber of this extension, see \cite[Corollary 3.3.6]{CP1}. Thus,
$$
\cX'=\cX_{\alpha,\alpha_2,\ldots,\alpha_m}=\cX_{\alpha,\beta_2,\ldots,\beta_m},
$$
which completes the proof of the claim.\medskip

Now from formula \eqref{210520j} and the fact that $\beta_i\in\Delta$ for all~$i$, it easily follows that
\begin{align*}
\cX'_\Delta = & (\cX_{\alpha,\alpha_2,\ldots,\alpha_m})_\Delta=(\cX_{\alpha,\beta_2,\ldots,\beta_m})_\Delta \\
= & ((\cX_{\alpha,\beta_2,\ldots,\beta_m})_\Delta )_{\beta_2,\ldots,\beta_m}\ge ((\cX_\alpha)_\Delta)_{\beta_2,\ldots,\beta_m}.
\end{align*}
The coherent configuration on the right-hand side of this relation is partly regular by Corollary~\ref{250720f}. Therefore, the coherent configuration $\cX'_\Delta$ is also partly regular. By  Theorem~\ref{200520i}, this implies that $\cX'_\Delta$ is schurian and separable, as required.\eprf

\section{A  lower bound for  indistinguishing number}\label{lb}

The main result of this section  (Theorem~\ref{201444a} below) establishes a lower bound for the indistinguishing number of a coherent configuration which is not partly regular, cf.~\cite[Theorem~3.1]{PV2017}. This bound gives a sufficient condition for a coherent configuration to be partly regular, and is used to prove Theorems~\ref{130520x} and~\ref{010620x} in the next section. 

\thrml{201444a}
Let $\cX$ be a coherent configuration of degree $n$, $k$ the maximal cardinality of a fiber of~$\cX$, and~$c=c(\cX)$.  If $\cX$ is not partly regular, then 
\qtnl{150620a}
(2k-1)c\ge n.
\eqtn
\ethrm 

\proof Let $\cX=(\Omega,S)$. In the sequel, $\Delta\in F(\cX)$ and $|\Delta|=k$. The fiber~$\Delta$ contains at least two points, for otherwise $k=1$ and $\cX$ is the discrete and hence partly regular configuration  in contrast to the hypothesis of the theorem. Set 
$$
\Delta_\alpha=\{\delta\in \Delta:\ n_{r(\alpha,\delta)}=1\}\qaq\Omega_1 =\{\alpha \in \Omega:\ \Delta_\alpha=\Delta\}.
$$

\lmml{250520a} 
$\Omega_1$ is a (possibly empty) union of fibers of cardinality~$k$. Moreover, the coherent configuration $\cX_{\Omega_1}$ is semiregular.
\elmm
\proof Let $\Gamma$  be a fiber containing a point of~$\Omega_1$.  The set $S_{\Gamma,\Delta}$ consists of~$s=r(\gamma,\delta)$, where $\gamma\in \Gamma\cap\Omega_1$ and $\delta$ runs over $\Delta$. Consequently, $n_s=1$ for all $s\in S_{\Gamma,\Delta}$. Therefore, $\Gamma\subseteq\Omega_1$. Thus, $\Omega_1$ is a union of fibers of $\cX$. Furthermore,  given $s\in S_{\Gamma,\Delta}$ we have
$$
k\ge |\Gamma|=n_{s^{}}\,|\Gamma|=n_{s^*}\,|\Delta|\ge |\Delta|=k.
$$ 
This proves the first statement. To prove  the second one, let $\Delta'$ and~$\Delta''$ be fibers contained in $\Omega_1$. By the definition of $\Omega_1$ and the first statement, any relations $s'\in S_{\Delta',\Delta}$ and $s''\in S_{\Delta,\Delta''}$ are matchings.  Therefore, $s'\cdot s''$ is a matching contained in $S_{\Delta',\Delta''}$. Thus, $S_{\Omega_1}$ consists of matchings and we are done.\eprf\medskip

By Lemma~\ref{250520a} and the hypothesis of the theorem,  the complement $\Omega'$  of the set~$\Omega_1$ contains at least two distinct points. 

\lmml{250520d}
For each $\gamma\in\Omega'$, 
$$
\sum_{s\in S_\gamma}n_s\ge \frac{k}{2},
$$
where $S_\gamma=\{r(\gamma,\delta):\ \delta\in \Delta$\, and\ \, $n_{r(\gamma,\delta)}>1\}$.
\elmm
\proof We have
\qtnl{160620a}
\sum_{s\in S_\gamma}n_s=\sum_{s\in S_\gamma}|\gamma s|=|\Delta|-|\Delta_\gamma|.
\eqtn
Since $|\Delta|=k$, this proves the required inequality if $\Delta_\gamma=\varnothing$. Let $\delta\in\Delta_\gamma$. It is easily seen that  $r(\delta,\lambda)=r(\delta,\gamma)\cdot r(\gamma,\lambda)$ is a matching of $\cX_\Delta$ for each $\lambda\in\Delta_\gamma$. Therefore,
\qtnl{260720a}
\Delta_\gamma\subseteq\{\lambda\in\Delta:\ n_{r(\delta,\lambda)}=1\}.
\eqtn
On the other hand, denote by $e$  the union of all matchings of the scheme $\cX_\Delta$. Then $e$ is a relation of this scheme. Moreover, $e$ is an equivalence relation on $\Delta$ (see~\cite[Theorem~2.1.25(4)]{CP1}) and the set on the right-hand side of~\eqref{260720a} is a class of~$e$.  In view of \cite[Corollary~2.1.23]{CP1}, the cardinality $a$ of this class divides $|\Delta|=k$. Furthermore, $a\ne k$, for otherwise, $\Delta_\gamma=\Delta$ and then $\gamma\in \Omega_1$, a contradiction. Thus,  
$$
|\Delta_\gamma|\le a\le \frac{k}{2}
$$ 
and the required statement follows from~\eqref{160620a}.\eprf

\lmml{220520o1}
Let $\varepsilon=\frac{2(k-1)}{2k-1}$. Assume that 
$
|\Omega'|\ge \varepsilon n.
$ 
Then inequality \eqref{150620a} holds.
\elmm
\proof  Denote by $N$ the cardinality of the set
\qtnl{130620v}
\{(\alpha,\beta,\gamma)\in\Delta\times\Delta\times\Omega':\ \alpha\ne\beta,\ \gamma\in c(\alpha,\beta)\},
\eqtn
where $c(\alpha,\beta)$ is as in formula~\eqref{100814c}. The number of $(\alpha,\beta)\in\Delta\times\Delta$ with $\alpha\ne\beta$ is equal to $k(k-1)$. Therefore there exists at least one such pair for which
\qtnl{130620j}
|c(\alpha,\beta)|\ge\frac{N}{k(k-1)}.
\eqtn

On the other hand, let $\gamma\in\Omega'$, and let $S_\gamma$ be as in Lemma~\ref{250520d}. Then $n_s\ge 2$ for all $s\in S_\gamma$. For  every such $s$ there are exactly $n_s(n_s-1)$ triples $(\alpha,\beta,\gamma)$ with distinct $\alpha,\beta\in \gamma s$, and all these triples belong to the set~\eqref{130620v}. By Lemma~\ref{250520d} this implies that
$$
N
=\sum_{\gamma\in\Omega'}\sum_{s\in S_{\Gamma,\Delta}}n_s(n_s-1)
\ge \sum_{\gamma\in\Omega'}\sum_{s\in S_\gamma}n_s
\ge \sum_{\gamma\in\Omega'}\frac{k}{2}= \frac{|\Omega'|\,k}{2},
$$
where $\Gamma$ is the fiber containing~$\gamma$. By formula~\eqref{130620j}  and the lemma assumption, we obtain
$$
c\ge |c(\alpha,\beta)|\ge \frac{N}{k(k-1)}\ge\frac{|\Omega'|}{2(k-1)}\ge
\frac{2(k-1)n}{2k-1}\cdot \frac{1}{2(k-1)}=\frac{n}{2k-1},
$$
as required.\eprf\medskip

By Lemma~\ref{220520o1}, we may assume that $|\Omega'|<\varepsilon n$. The coherent configuration $\cX$ is not partly regular. Therefore no point $\delta\in\Omega_1$ is  regular and there exist distinct  $\alpha,\beta\in\Omega'$ such that $\delta\in c(\alpha,\beta)$. Since the coherent configuration $\cX_{\Omega_1}$ is semiregular (Lemma~\ref{250520a}),   the relation $s=r(\delta,\lambda)$ is a matching for all $\lambda\in\Omega_1$. It follows that
$$
r(\alpha,\lambda)=r(\alpha,\delta)\cdot s=r(\beta,\delta)\cdot s=r(\beta,\lambda).
$$
Consequently, $\Omega_1\subseteq c(\alpha,\beta)$. This implies that
$$
c\ge |c(\alpha,\beta)|\ge |\Omega_1|=n-\varepsilon n>n\left(1-\frac{2(k-1)}{2k-1}\right)=\frac{n}{2k-1}
$$
which completes the proof.\eprf

\crllrl{423939b}
Let $\cX$ be a coherent configuration of degree $n$, $c=c(\cX)$, and  $t$  an irreflexive basis relation of~$\cX$. Assume that $(2m_t-1)\,c< n$, where
\qtnl{280520o}
m_t=\max\limits_{r,s \in S} c_{r s}^t.
\eqtn
Then the extension of $\cX$ with respect to any two points  forming a pair from~$t$ is partly regular. 
\ecrllr

\proof Let $\cX'$ be the extension of $\cX$ with respect to the points $\alpha,\beta$ such that $(\alpha, \beta) \in t$.  Then each fiber $\Delta$ of $\cX'$ different from both $\{\alpha\}$ and $\{\beta\}$ is contained in the set $\alpha r \cap \beta s^*$ for appropriate $r, s \in S$. It follows that
$$
|\Delta|\le  |\alpha r \cap \beta s^*|=c_{r s}^t \leq m_t. 
$$
Thus the maximal cardinality $k'$ of a fiber of $\cX'$ is less than or equal to $m_t$. Since obviously $c'=c(\cX')$ is less than or equal to~$c$, the condition of the corollary implies that
$$
(2k'-1)\,c'\le  (2m_t-1)\,c< n.
$$ 
Thus $\cX'$  is  partly regular by Theorem~\ref{201444a}.\eprf

\section{Small Hollmann  and  Passman schemes}\label{sH}

\subsection{Algebraic fusion.}
Let $\cX=(\Omega, S)$ be a coherent configuration, and  let $\Phi$ be a group of algebraic automorphisms of~$\cX$. For each $s \in S$, set
$$
s^\Phi=\bigcup_{\varphi \in \Phi}\varphi(s).
$$
Clearly, $(1_\Omega)^\Phi=1_\Omega$. Moreover the set $S^\Phi=\{s^\Phi:\ s\in S\}$ forms a partition of the Cartesian square $\Omega^2$. According to \cite[Lemma~2.3.26]{CP1},  the pair 
$$
\cX^\Phi=(\Om, S^\Phi)
$$ 
is a coherent configuration called the algebraic fusion of $\cX$ with respect to $\Phi$. In the following lemma, we establish a simple upper bound for the intersection numbers of an algebraic fusion.

\lmml{411958b}
In the above notation, let $r, s, t\in S$ and $m_t$ is as in~\eqref{280520o}. Then
$$
c_{r^{\Phi} s^{\Phi}}^{t^{\Phi}} \leq m_t\,|\Phi|^{2}.
$$
\elmm
\proof We have 
$c_{r^{\Phi} s^{\Phi}}^{t^{\Phi}} 
 =\sum\limits_{\varphi, \psi \in \Phi} c_{\varphi(r)\, \psi(s)}^t \leq m_t\,|\Phi|^{2}.
 $
\eprf

\subsection{Proof of Theorem~\ref{130520x}.}
Let $\cX$ be a small Hollmann scheme. Then the degree of $\cX$ is equal to $n=q(q-1)/2$, where $q = 2^d$ for a prime $d\ge 3$. Moreover,~$\cX$ is associated with the permutation group $G=\PGaL(2,q)$ of degree~$n$ from  Theorem~\ref{140520c}(1). As in Subsection~\ref{260720x}, one can see that $\cX$ coincides with symmetric pseudocyclic scheme $\cX'$ of degree $n$ and valency  $d(q+1)$, associated with the group~$\PGaL(2,q)$ and studied~\cite{HX2}. In particular, $\cX$ is obtained from the large Hollmann scheme of degree $n$ by merging the basis relations via the Frobenius map $x\mapsto x^2$, $x\in\mF_q$. In other words, $\cX$ is the algebraic fusion of the large Hollmann scheme of degree $n$ with respect the induced action of $\aut(\mF_q)$ on its basis relations.

\prpstnl{270520i}
Let $\cX$ and $\cX_q$ be the small and large Hollmann schemes of  degree $q(q-1)/2$, respectively. Then 
\nmrt
\tm{1} $\cX={\cX_q}^{\hspace{-2pt}\Phi}$, where $\Phi\le\aut_{alg}(\cX_q)$ is a group of order $d$,
\tm{2} $\cX$  is a pseudocyclic scheme of valency~$d(q+1)$,
\tm{3} for each irreflexive $t\in S(\cX)$, we have $m_t\le 4d^2$, where $m_t$ is as in~\eqref{280520o}. 
\enmrt
\eprpstn
\proof Statements~(1) and (2) follow from the above discussion. 
Next,  from the formulas for the intersection numbers of the scheme $\cX_q$, given in~\cite[Theorem 2.2]{HX2},  it follows that $c_{xy}^z\le 4$ for all irreflexive $x,y,z\in S(\cX_q)$. In particular, 
$$
m_z\le 4. 
$$
On the other hand, by statement~(1), each irreflexive $t\in S(\cX)$ is of the form~$z^\Phi$ for some irreflexive~$z$. Thus by Lemma~\ref{411958b},
$$
m_t=\max\limits_{x,y \in S(\cX_q)} c_{x^\Phi y^\Phi}^{z^\Phi}\le m_z |\Phi|^2\le 4d^2,
$$
which  proves statement~(3).\eprf\medskip

Let us prove Theorem~\ref{130520x}. If $d=3$, then a straightforward calculation shows that $\cX$ is trivial scheme and hence $c(\cX)=1$.  Let $d>3$. By Lemma~\ref{030620d} and Theorem~\ref{200520i}, it suffices to verify that the extension of  $\cX$ with respect to at least two  points is partly regular. 

\thrml{280520a}
Let $q=2^d$, where $d>3$ is a prime and $d\ne 7,11,13$.  Then every extension of the small Hollmann scheme of degree $q(q-1)/2$ with respect to at least two points is partly regular. 
\ethrm
\proof Let $\cX$ be the small Hollmann scheme of  degree $n=q(q-1)/2$. By Proposition~\ref{270520i}(2), the number $c=c(\cX)$ is equal to $d(2^d+1)-1$. By statement~(3) of the same proposition, $m_t\le 4\,d^2$ for any irreflexive $t\in S(\cX)$. Now if $d> 16$, then
$$
(2m_t-1)\,c< (8d^2-1)\,(d\,(2^d+1)-1) <2^{d-1}(2^d-1)=n.
$$
By Corollary~\ref{423939b}, this proves the required statement for all (prime) $d>13$. In the remaining case, $d=5$, the required statement has been checked with the help of the computer package~COCO2P~\cite{KlinCOCO2P}.\eprf

\subsection{Proof of Theorem~\ref{010620x}.}
Let $q$ be an odd prime power. The permutation group $G\le\AGL(2,q)$ defined in Theorem~\ref{140520c}(2) has a Frobenius subgroup $H$ consisting of the permutations
\qtnl{310520a}
\begin{pmatrix}
x\\ y\\
\end{pmatrix}
\rightarrow
\begin{pmatrix}
a & 0 \\ 0 & a^{-1}\\ 
\end{pmatrix}
\begin{pmatrix}
x\\ y\\
\end{pmatrix}
+
\begin{pmatrix}
b\\ c\\
\end{pmatrix},\qquad  x,y\in \mF_q,
\eqtn
where $a, b, c \in \mF_q$ and $a \ne 0$. The group $D\le\GL(2,q)$ consisting of the four matrices
$$
\begin{pmatrix}
\pm 1 & 0 \\ 0 & \pm 1\\ 
\end{pmatrix}
\qaq
\begin{pmatrix}
0 & \pm 1 \\ \pm 1 & 0\\ 
\end{pmatrix}
$$
is contained in $G$, normalizes $H$, and, moreover, $G=H\rtimes D$. In particular, $D$ acts on the basis relations of the Frobenius scheme $\cY=\inv(H)$ as a group $\Phi$ of algebraic automorphisms.

\prpstnl{300520a}
Let $\cX=\inv(G)$ be the Passman scheme of  degree $q^2$. Then 
\nmrt
\tm{1} $\cX=\cY^{\Phi}$, where $\Phi\le\aut_{alg}(\cY)$ is a group of order $2$,
\tm{2} $\cX$  is a pseudocyclic scheme of valency~$2(q-1)$,
\tm{3} there is  irreflexive $t\in S(\cX)$ such that $m_t\le 4$, where $m_t$ is as in~\eqref{280520o}. 
\enmrt
\eprpstn
\proof Statements~(1)  and (2) follow from \cite[Section~4.3]{MP1}. To prove statement~(3), let $u$ be the basis relation of~$\cY$, containing the pair $(\alpha,\beta)$, where 
$$
\alpha=(0,0)\qaq \beta=(1,1).
$$ 
It suffices to verify that $m_u=1$; indeed, then $t=u^\Phi$ is the required relation by statement~(1) and Lemma~\ref{411958b}.\medskip

We need to verify that $c_{r s}^u \leq 1$  for all $r, s \in S(\cY)$.  Without loss of generality, we may assume that  $r$ and $s$ are such that $c_{r s}^u \ne 0$. Then there exists $\gamma\in \alpha r\cap\beta s^*$ for which $\alpha r=\gamma^{H_\alpha}$, $\beta s^*=\gamma^{H_\beta}$, and 
\qtnl{260720z}
c_{r s}^u=|\gamma^{H_\alpha} \cap \gamma^{H_\beta}|.
\eqtn
Let us calculate the number on the right-hand side. Using the explicit form~\eqref{310520a} of the elements of~$H$, one can easily find that the groups $H_\alpha$ and $H_\beta$ consist of permutations the parameters $a,b,c$ of which satisfy the relations
$$
b=c=0 \qaq a+b=a^{-1}+c=1,
$$
respectively. Consequently, assuming $\gamma=(x,y)$, we have
$$
\gamma^{H_\alpha}=\left\{\left(ax, \frac{y}{a}\right):\  a \in \mF_q^*\right\}
\qaq
\gamma^{H_\beta}=\left\{\left(ax +1-a, \frac{y}{a}+1-\frac{1}{a}\right): a \in \mF_q^*\right\}.
$$
In view of~\eqref{260720z}, the intersection number $c_{r s}^u$ is equal to the number of elements $a \in \mF_q^*$ such that
$$
ax=ax+1-a \qaq \frac{y}{a}=\frac{y}{a}+1-\frac{1}{a}.
$$
These equations are satisfied for $a=1$ only. Therefore, $c_{rs}^u=1$, as required.\eprf\medskip

Let us prove Theorem~\ref{010620x}. By Lemma~\ref{030620d} and Theorem~\ref{200520i}, it suffices to verify that the extension of  $\cX$ with respect to at least two  points is partly regular. 

\thrml{310520d}
Let $q$ be an odd prime power. Then every extension of the Passman scheme of degree $q^2$ with respect to at least two points is partly regular. 
\ethrm
\proof Let $\cX$ be the Passman scheme of  degree $n=q^2$. By Proposition~\ref{300520a}(2), the number $c=c(\cX)$ is equal to $2(q-1)-1$. Let $t$ be the basis relation of $\cX$, defined in  Proposition~\ref{300520a}(3). Then $m_t\le 4$. Now if $q\ge 13$, then
$$
(2m_t-1)\,c \le 7\, (2q-3) <q^2=n.
$$
By Corollary~\ref{423939b}, this proves the required statement for all $q\ge 13$. In the remaining case, the required statement has been checked with the help of the computer package~COCO2P~\cite{KlinCOCO2P}.\eprf

\section{Concluding remarks and open problems}\label{030620a}

\subsection{Pseudocyclic schemes.}\label{110620a} A scheme is said to be {\it $k$-equivalenced} (and just equivalenced if~$k$ is irrelevant) if all irreflexive basis relations of it have  valency~$k$. It is known that every $k$-equivalenced scheme is pseudocyclic for $1\le k\le 4$; this follows from  results obtained in~\cite{MZ,Park2015,Moon20} and~\cite[Theorem 3.1]{MP1}. \medskip

By Theorem~\ref{140520c}, if $\cX$ is a schurian equivalenced scheme of sufficiently large degree, which is not  trivial or Frobenius,  then  either $\cX$ is exceptional, i.e., the Hollmann or Passman scheme,  or the inclusion $\aut(\cX)\le\AGaL(1,q)$ holds for some~$q$. In all   cases except for the last one, $\cX$ is pseudocyclic, see Propositions~\ref{160520i}, \ref{270520i}(2), and~\ref{300520a}(2), and~\cite[Theorem~3.1]{MP1}. In the latter case, the group $\AGaL(1,q)$ can contain  $\frac{3}{2}$-transitive subgroups which are not $2$-equivalent to Frobenius groups (such a subgroup is always primitive). We do not know whether the scheme of at least one of these subgroups is not pseudocyclic.

\subsection{Separability number.} Finding the exact values of $s(\cX)$ for an exceptional scheme $\cX$ is still an open problem. A direct calculation shows that these schemes are separable for small~$q$.

\subsection{Superschurian schemes.} The following concept was first formulated many years ago in discussions of the third author with Sergei Evdokimov. A scheme $\cX$ is said to be {\it superschurian} if the extension of $\cX$ with respect to every set of points is schurian. In particular, all  superschurian schemes are schurian.  In fact, only a few families of superschurian schemes are known; these include partly regular schemes, cyclotomic schemes over finite fields, normal circulant schemes, and some TI-schemes \cite{EP1,MP1,CP2}. Theorem~\ref{030620i} implies that any large Hollmann scheme is superschurian. We do not know whether other exceptional schemes are superschurian.

\subsection{Base number.} The base number $b(\cX)$ of a coherent configuration $\cX$ is defined to be the smallest number of points such that the extension of $\cX$ with respect to them is the discrete configuration, see \cite[Section~3.3.2]{CP1}. In general,  the base number of the group $\aut(\cX)$ is less than or equal to $b(\cX)$. The equality is attained, for example, if $\cX$ is a partly regular coherent configuration. By virtue of this observation,  Theorems~\ref{030620i}, \ref{280520a}, and \ref{310520d} imply that except, possibly, for several small Hollmann schemes the equality holds also for all exceptional pseudocyclic schemes.\medskip

In fact, the base number of an exceptional scheme of enough large degree is bounded by $3$. This fact can be used to construct a polynomial-time algorithm recognizing whether or not a given scheme is exceptional. Taking the above discussion into account, this reduces  the recognition problem for the class of schurian equivalenced schemes (see \cite[p.281]{VC}) to the non-Frobenius schemes $\cX$ for which  $\aut(\cX)\le\AGaL(1,q)$.

\subsection{Bound in Theorem~\ref{201444a}.}  Denote by  $f(n)$ the maximum of the ratio $\frac{n}{c(\cX)k(\cX)}$ taken over all  non-partly-regular coherent configurations $\cX$ of degree $n$, where $k(\cX)$ is the maximal cardinality of a fiber of~$\cX$. Clearly, $f(n)>0$. Theorem~\ref{201444a} states that $f(n)<2$. It would be interesting to find the function $f(n)$ explicitly. We have found a (schurian) non-partly-regular coherent configuration with parameters
$$
n=24,\quad k=8,\quad c=4,
$$
which shows that $f(24)\ge 3/4$.

\end{document}